\documentclass[11pt]{article}
\usepackage{amssymb,amsmath,amsthm, amsfonts}
\usepackage{graphicx}
\usepackage{epsfig}
\usepackage{color}

\def\dref#1{(\ref{#1})}

\usepackage{amssymb,amsmath,amsthm, amsfonts, mathrsfs, bbm, dsfont}

\theoremstyle{plain}

\theoremstyle{definition}

\numberwithin{equation}{section}

\begin{document}

\title{{Influence of flux limitation on large time behavior in a three-dimensional chemotaxis-Stokes system modeling coral fertilization}}

\author{
{\rm Ji Liu$^{*}$
}\\[0.2cm]
{\it\small \rm College of Sciences, Nanjing Agricultural University}\\
{\it\small \rm Nanjing 210095, People's Republic of China}\\
}
\date{}

\maketitle

\renewcommand{\thefootnote}{\fnsymbol{footnote}}
\setcounter{footnote}{-1}
\footnote{$^{*}$Corresponding author.
\\E-mail addresses: Liuji@njau.edu.cn(J. Liu).}

\maketitle \vspace{0.3cm} \noindent

\begin{minipage}[t]{15.5cm}
{\bf Abstract:} In this paper, we consider the following system
$$\left\{\begin{array}{ll}
n_t+u\cdot\nabla n&=\Delta n-\nabla\cdot(n\mathcal{S}(|\nabla c|^2)\nabla c)-nm,\\
c_t+u\cdot\nabla c&=\Delta c-c+m,\\
m_t+u\cdot\nabla m&=\Delta m-mn,\\
u_t&=\Delta u+\nabla P+(n+m)\nabla\Phi,\qquad \nabla\cdot u=0
\end{array}\right.$$
which models the process of coral fertilization, in a smoothly three-dimensional bounded domain, where $\mathcal{S}$ is a given function fulfilling
$$|\mathcal{S}(\sigma)|\leq K_{\mathcal{S}}(1+\sigma)^{-\frac{\theta}{2}},\qquad \sigma\geq 0$$
with some $K_{\mathcal{S}}>0.$ Based on conditional estimates of the quantity $c$ and the gradients thereof, a relatively compressed argument as compared to that proceeding in related precedents shows that if
$$\theta>0,$$
then for any initial data with proper regularity an associated initial-boundary problem under no-flux/no-flux/no-flux/Dirichlet boundary conditions admits a unique classical solution which is globally bounded, and which also enjoys the stabilization features in the sense that
$$\|n(\cdot,t)-n_{\infty}\|_{L^{\infty}(\Omega)}+\|c(\cdot,t)-m_{\infty}\|_{W^{1,\infty}(\Omega)}
+\|m(\cdot,t)-m_{\infty}\|_{W^{1,\infty}(\Omega)}+\|u(\cdot,t)\|_{L^{\infty}(\Omega)}\rightarrow0
\quad\textrm{as}~t\rightarrow \infty$$
with $n_{\infty}:=\frac{1}{|\Omega|}\left\{\int_{\Omega}n_0-\int_{\Omega}m_0\right\}_{+}$ and $m_{\infty}:=\frac{1}{|\Omega|}\left\{\int_{\Omega}m_0-\int_{\Omega}n_0\right\}_{+}.$
\end{minipage}

\vspace{0.3cm}
\vspace{0.3cm}






\section{Introduction}

As shown in the experiments \cite{Coll1,Coll2,Miller1,Miller2}, during the process of coral fertilization, chemotaxis phenomenon may occur in terms of oriented motion of sperms in response to some kind of chemical signal released by eggs. In mathematics, it can be modeled by
\begin{equation}
\left\{\begin{array}{ll}
n_t+u\cdot\nabla n&=\Delta n-\nabla\cdot(n\mathcal{S}(n,c)\nabla c)-nm,\\
c_t+u\cdot\nabla c&=\Delta c-c+m,\\
m_t+u\cdot\nabla m&=\Delta m-mn,\\
u_t+\kappa(u\cdot\nabla)u&=\Delta u+\nabla P+(n+m)\nabla\Phi,
\end{array}\right.\label{1.1}
\end{equation}
where both $n$ and $m$ denote the densities of unfertilized sperms and eggs, respectively, $c$ represents the concentration of the signal, $u$ stands for the velocity of the ambient ocean flow, $P$ expresses the pressure within the fluid and $\Phi$ is a given function used to denote the gravitational potential (\cite{Espejo1}).

In the simplified case when $n\equiv m,$ mathematical results on global dynamics of system \dref{1.1} can be found in \cite{Kiselev 1,Kiselev 2,Espejo}, where, in particular, it is shown in \cite{Kiselev 1,Kiselev 2} that with given fluid velocity the thoroughness level of the fertilization depends on whether or not the role of chemotaxis comes into play, and where if the fluid flow is slow, i.e. $\kappa=0$ in \dref{1.1}, but the velocity thereof is unknown, then system \dref{1.1} is globally solvable in the weak sense (\cite{Espejo}). Whereas for the more realistic situation that $n\not\equiv m,$ complete consumption of eggs needs not only sufficiently high concentration of signal but also adequately large initial densities for both sperms and eggs (\cite{Espejo2}).

Thereafter, the global dynamics of system \dref{1.1} with $n\not\equiv m$ and with an unknown fluid velocity $u,$ which is referred as a more challenging case, becomes an interesting subject. Specially, in spatially two-dimensional setting, it is proved that system \dref{1.1} with $\kappa=1$ possesses a unique classical solution which is globally bounded and approaches to some constant equilibrium as time goes to infinity (\cite{Espejo1}). While for physically most relevant three-dimensional case, the corresponding global solvability of \dref{1.1} in the context of Stokes-fluid or Navier--Stokes-fluid requires some smallness hypothesis for initial data (\cite{Li2,Li3,Chae,Htwe}), or necessary aids of some type of nonlinear mechanism, such as porous medium diffusion (\cite{Liu3}), or signal-dependent sensitivity (\cite{Li1}), or saturation effects of cells (\cite{Li2,Htwe,Liu0,Zheng}), or $p$-Laplace diffusion of cells (\cite{Liu2}).

Recently, from some refined models proposed in \cite{Bellomo,Bianchi,Bianchi1}, it can be observed that the migration of cells rests with some gradient-dependent limitations, which inspires investigations of global dynamics on associated initial-boundary problems. In particular, it is uncovered that global solvability can be enforced by appropriately strong p-Laplace type cell diffusion (\cite{Bendahmane,LiY,TaoW,TaoW1,Zhuang,Liu2}), and that the precise parameters which describe the critical interaction of flux-limited cross diffusion and relativistic heat equation type diffusion can be achieved by constructions of suitable sub-solutions on the basis of a comparison argument (\cite{Bellomo1,Bellomo2,Mizukami,Chiyoda}).

When merely flux limitation of signal is taken into account, that is the function $\mathcal{S}$ in \dref{1.1} is supposed to generalize the prototype
\begin{equation}
\mathcal{S}(\sigma)=K_{\mathcal{S}}(1+\sigma)^{-\frac{\theta}{2}},\quad\textrm{for any}~\sigma\geq 0 \label{1.2}
\end{equation}
with $K_{\mathcal{S}}>0$ and $\theta>0,$ the simplified parabolic-elliptic Keller--Segel system under homogeneous Neumann boundary conditions is globally solvable with a bounded solution in the classical sense under the assumptions that
\begin{equation}
|\mathcal{S}(\sigma)|\leq K_{\mathcal{S}}(1+\sigma)^{-\frac{\theta}{2}},\quad\textrm{for any}~\sigma\geq 0 \label{1.3}
\end{equation}
with $\theta=\theta_c:=\frac{N-2}{N-1},$ where $N\geq2$ denotes the spatial dimension. Whereas on the contrary, if
$$\mathcal{S}(\sigma)>k_{\mathcal{S}}(1+\sigma)^{-\frac{\theta}{2}},\quad\textrm{for any}~\sigma\geq 0$$
with some $k_{\mathcal{S}}>0$ and $\theta<\theta_c,$ spontaneous aggregation in the sense of unbounded density may emerge at some finite time for $N\geq 3$ (\cite{WinklerA}). For more complex case that the liquid circumstance is extra accounted for (\cite{WinklerC}), the admissible range of $\theta$ for global solvability coincides with \dref{1.3} in spatially three-dimensional setting, i.e.
$$\theta>\frac{1}{2}.$$

In view of the above results, we intend to discover that to which extent the saturation effects induced by $\mathcal{S}$ fulfilling \dref{1.3} can prevent the phenomenon of explosion to occur spontaneously to an associated initial-boundary problem with the coral fertilization model \dref{1.1}. Particularly, the precise problem we consider herein is
\begin{equation}
\left\{\begin{array}{ll}
n_t+u\cdot\nabla n=\Delta n-\nabla\cdot(n\mathcal{S}(|\nabla c|^2)\nabla c)-nm,&x\in\Omega,~t>0,\\
c_t+u\cdot\nabla c=\Delta c-c+m,&x\in\Omega,~t>0,\\
m_t+u\cdot\nabla m=\Delta m-mn,&x\in\Omega,~t>0,\\
u_t=\Delta u+\nabla P+(n+m)\nabla\Phi,\qquad \nabla\cdot u=0,&x\in\Omega,~t>0,\\
\frac{\partial n}{\partial \nu}=\frac{\partial c}{\partial \nu}=\frac{\partial m}{\partial \nu}=0,\qquad u=0,&x\in\partial\Omega,~t>0,\\
n(x,0)=n_0(x),~c(x,0)=c_0(x),~m(x,0)=m_0(x),~u(x,0)=u_0(x),&x\in\Omega,
\end{array}\right.\label{1.4}\end{equation}
where $\Omega\subset\mathbb{R}^3$ is a general bounded domain with smooth boundary, where the gravitational potential $\Phi$ complies with
\begin{equation}
\Phi\in W^{2,\infty}(\Omega),\label{1.5}
\end{equation}
and where for certain $K_0>0$ the initial data $(n_0,c_0,m_0,u_0)$ fulfills
\begin{equation}\left\{\begin{array}{ll}
n_0\in C^0(\bar\Omega)~~\textrm{such that}~~n_0\geq 0~~\textrm{and}~~\|n_0\|_{C^0(\Omega)}\leq K_0,\\
c_0\in W^{1,\infty}(\Omega)~~\textrm{such that}~~c_0\geq 0~~\textrm{and}~~\|c_0\|_{W^{1,\infty}(\Omega)}\leq K_0,\\
m_0\in W^{1,\infty}(\Omega)~~\textrm{such that}~~m_0\geq 0~~\textrm{and}~~\|m_0\|_{W^{1,\infty}(\Omega)}\leq K_0,\\
u_0\in \bigcup_{\alpha\in(\frac{3}{4},1)}D(A^{\alpha})~~\textrm{such that}~~\|A^{\alpha}u_0\|_{L^2(\Omega)}\leq K_0
\end{array}\right.\label{1.6}\end{equation}
with $A:=-\mathcal{P}\Delta$ standing for the realization of the Stokes operator whose domain is $D(A):=W^{2,2}(\Omega;\mathbb{R}^3)\bigcap W^{1,2}_{0,\sigma}(\Omega),$ where $\mathcal{P}$ denotes the Helmholtz projection on $L^{2}(\Omega;\mathbb{R}^3),$ and where $W^{1,2}_{0,\sigma}(\Omega):=W^{1,2}_{0}(\Omega;\mathbb{R}^3)\bigcap L^2_{\sigma}(\Omega)$ with $L^2_{\sigma}(\Omega):=\{\omega\in L^{2}(\Omega;\mathbb{R}^3)|\nabla\cdot\omega=0\}.$

Based on well-established conditional regularity features of the fluid field $u$ (\cite{WinklerC}), a suitable application of smoothness estimates for Neumann heat semigroup can establish an inequality which shows a relationship between uniform $L^{\infty}$ estimates of the signal gradient and the conditional $L^p$ bounds of $n,$ which underlies the derivation of the following results on global boundedness (in Section 5).

\textbf{Theorem 1.1} \emph{Assume that both \dref{1.5} and \dref{1.6} are valid. Let $\mathcal{S}\in C^2([0,\infty))$ fulfill
\begin{equation}
|\mathcal{S}(\sigma)|\leq K_{\mathcal{S}}(1+\sigma)^{-\frac{\theta}{2}},\quad\textrm{for any}~\sigma\geq 0
\label{1.7}\end{equation}
with certain $K_{\mathcal{S}}>0$ and
$$\theta>0.$$
Then one can find a unique quadruple of functions
\begin{equation}\left\{\begin{array}{ll}
n\in C^0(\bar\Omega\times[0,\infty))\bigcap C^{2,1}(\bar\Omega\times(0,\infty)),\\
c\in \bigcap_{l>3} C^0([0,\infty);W^{1,l}(\Omega))\bigcap C^{2,1}(\bar\Omega\times(0,\infty)),\\
m\in \bigcap_{r>3} C^0([0,\infty);W^{1,r}(\Omega))\bigcap C^{2,1}(\bar\Omega\times(0,\infty)),\\
u\in \bigcup_{\alpha\in(\frac{3}{4},1)}C^0([0,\infty);D(A^{\alpha}))\bigcap C^{2,1}(\bar\Omega\times(0,\infty);\mathbb{R}^3)
\end{array}\right.
\label{1.8}\end{equation}
with $n\geq0,c\geq0$ and $m\geq0$ in $\Omega\times(0,\infty)$ such that combined with some $P\in C^{1,0}\left(\Omega\times(0,\infty)\right),$ the quintuple $(n,c,m,u,P)$ constitutes a classical solution of \dref{1.4} in $\Omega\times(0,\infty),$ and is also bounded in line with
\begin{equation}\|n(\cdot,t)\|_{L^{\infty}(\Omega)}+\|c(\cdot,t)\|_{W^{1,\infty}(\Omega)}
+\|m(\cdot,t)\|_{W^{1,\infty}(\Omega)}+\|A^{\alpha}u(\cdot,t)\|_{L^2(\Omega)}\leq C\quad\textrm{for all}~t>0
\label{1.9}\end{equation}
with some $\alpha\in(\frac{3}{4},1)$ and $C>0.$}

Relying on the boundedness features in various forms with regards to each component of the solution and on the basic relaxation properties exhibited in Lemma 5.1 below, we are also able to achieve the following statement on the convergence for each component of the solution as time becomes arbitrarily large.

\textbf{Theorem 1.2} \emph{The global classical solution constructed in Theorem 1.1 has the property that
\begin{equation}\|n(\cdot,t)-n_{\infty}\|_{L^{\infty}(\Omega)}+\|c(\cdot,t)-m_{\infty}\|_{W^{1,\infty}(\Omega)}
+\|m(\cdot,t)-m_{\infty}\|_{W^{1,\infty}(\Omega)}+\|u(\cdot,t)\|_{L^{\infty}(\Omega)}\rightarrow0
\quad\textrm{as}~t\rightarrow\infty,
\label{1.10}\end{equation}
where $n_{\infty}:=\frac{1}{|\Omega|}\left\{\int_{\Omega}n_0-\int_{\Omega}m_0\right\}_{+}$ and $m_{\infty}:=\frac{1}{|\Omega|}\left\{\int_{\Omega}m_0-\int_{\Omega}n_0\right\}_{+}.$}

More precisely, with the aids of the basic relaxation features in Lemma 5.1 as well as the uniform $L^{\infty}$ bounds of $n$ provided by Theorem 1.1, some constant equilibrium can be detected at first for the solution which converges in space with lower regularity, such as $L^2$ space (see subsection 5.1). In subsequence, by means of an associated Ehrling lemma, this in conjunction with the H\"{o}lder regularities of both $n$ and $u$ as well as of the gradients of both $m$ and $c$ improves the regularity of the space, in which the solution stabilizes towards the constant equilibrium, to the desired level as asserted in Theorem 1.2 (see subsection 5.2).

\section{Preliminaries}

In light of the treatments for closely related problems (\cite{Winkler 3}), the following assertions on local solvability and extensibility of \dref{1.4} are valid, though accounted for the evolution of one more quantity $m.$

\textbf{Lemma 2.1} \emph{Let the initial data $(n_0,c_0,m_0,u_0)$ comply with \dref{1.6}, and let both \dref{1.5} and \dref{1.7} hold. Then one can find $T_{\max}\in(0,+\infty]$ and a unique quadruple of functions $(n,c,m,u)$ satisfying
\begin{equation}\left\{\begin{array}{ll}
n \in C^0(\bar\Omega\times[0,T_{\max}))\bigcap C^{2,1}(\bar\Omega\times(0,T_{\max})),\\
c \in\bigcap_{l>3}C^0([0,T_{\max});W^{1,l}(\Omega))\bigcap C^{2,1}(\bar\Omega\times(0,T_{\max})),\\
m \in\bigcap_{r>3}C^0([0,T_{\max});W^{1,r}(\Omega))\bigcap C^{2,1}(\bar\Omega\times(0,T_{\max})),\\
u \in \bigcup_{\alpha\in(\frac{3}{4},1)}C^0([0,T_{\max});D(A^{\alpha}))\bigcap C^{2,1}(\bar\Omega\times(0,T_{\max});\mathbb{R}^3),
\end{array}\right.\end{equation}\label{2.1}
and $n\geq0,~c\geq0$ in $\Omega\times(0,T_{\max}),$ such that there exists some $P\in C^{2,1}(\Omega\times(0,T_{\max}))$ which along with $(n,c,m,u)$ solves \dref{1.4} classically in $\Omega\times(0,T_{\max}),$ and that
\begin{equation}\begin{aligned}
&\textrm{either}\quad T_{\max}=\infty\quad\textrm{or}\quad\textrm{for any}\quad\alpha\in(\frac{3}{4},1),\\
&\lim_{t\nearrow T_{\max}}\sup\left\{\|n(\cdot,t)\|_{L^{\infty}(\Omega)}+\|c(\cdot,t)\|_{W^{1,l}(\Omega)}
+\|m(\cdot,t)\|_{W^{1,r}(\Omega)}+\|A^{\alpha}u(\cdot,t)\|_{L^2(\Omega)}\right\}
=\infty
\end{aligned}\label{2.2}\end{equation}
is valid.}

\textbf{Proof.} Lemma 2.1 follows from an appropriate modification of the reasoning of \cite[Lemma 2.1]{Winkler 3}.

Thanks to the nonnegativity of both $n$ and $m,$ we also have the basic estimates as follows.

\textbf{Lemma 2.2} \emph{Suppose $(n,c,m,u)$ is the solution as provided in Lemma 2.1. We have
\begin{equation}\int_{\Omega}n(\cdot,t)\leq \int_{\Omega}n_0,\quad\int_{\Omega}m(\cdot,t)\leq \int_{\Omega}m_0\quad\textrm{for each}\quad t\in(0,T_{\max})\label{2.3}\end{equation}
and
\begin{equation}\|m(\cdot,t)\|_{L^{\infty}(\Omega)}\leq \|m_0\|_{L^{\infty}(\Omega)}\quad\textrm{for each}\quad t\in(0,T_{\max})\label{2.4}\end{equation}
as well as
\begin{equation}\|c(\cdot,t)\|_{L^{\infty}(\Omega)}\leq \max\{\|c_0\|_{L^{\infty}(\Omega)},\|m_0\|_{L^{\infty}(\Omega)}\}\quad\textrm{for each}\quad t\in(0,T_{\max}).\label{2.5}\end{equation}}

\textbf{Proof.} In view of the nonnegativity of $n$ and $m,$ \dref{2.3} follows from straightforward integrations of $n$-equation and $m$-equation in \dref{1.4}, respectively. Also due to $n\geq 0$ and $m\geq 0,$ \dref{2.4} can be inferred from the maximum principle. By virtue of \dref{2.4}, \dref{2.5} is a consequence of the comparison principle.

The following conditional estimates of the fluid velocity, which are derived from well-established reasoning frameworks in \cite[Proposition 1.1 and Corollary 2.1]{WinklerC}, will play a helpful role in the analysis of Section 3.

\textbf{Lemma 2.3} \emph{Let $(n,c,m,u)$ be a solution constructed in Lemma 2.1. Then for some $\alpha\in(\frac{3}{4},1),p\geq 2,q>3$ and $\delta>0,$ there exist $K_1=K_1(\alpha,p,q,\delta)>0$ and $K_2=K_2(\alpha,p,q,\delta)>0$ with the properties that
\begin{equation}\|A^{\alpha}u(\cdot,t)\|_{L^2(\Omega)}\leq K_1\cdot\left\{1+\sup_{s\in(0,t)}\|n(\cdot,s)\|_{L^p(\Omega)}\right\}^{\frac{p}{p-1}\cdot\left(\frac{4\alpha-1}{6}+\delta\right)}\quad\textrm{for any}\quad t\in(0,T_{\max})\label{2.6}\end{equation}
and ,whereafter, that
\begin{equation}\|u(\cdot,t)\|_{L^q(\Omega)}\leq K_2\cdot\left\{1+\sup_{s\in(0,t)}\|n(\cdot,s)\|_{L^p(\Omega)}\right\}^{\frac{p}{p-1}\cdot\left(\frac{q-3}{3q}+\delta\right)}\quad\textrm{for any}\quad t\in(0,T_{\max}).\label{2.7}\end{equation}}

\textbf{Proof.} Based on $u$-equation in \dref{1.4}, it follows from the reasoning of \cite[Proposition 1.1]{WinklerC} that
\begin{equation}\|A^{\alpha}u(\cdot,t)\|_{L^2(\Omega)}\leq C_1\cdot\left\{1+\sup_{s\in(0,t)}\|n(\cdot,s)+m(\cdot,s)\|_{L^p(\Omega)}\right\}^{\frac{p}{p-1}\cdot\left(\frac{4\alpha-1}{6}+\delta\right)}
\label{2.8}\end{equation}
with some $C_1=C_1(K_0,|\Omega|,p,\alpha,\delta)>0$ for any $\in(0,T_{\max}).$ From \dref{1.6} and \dref{2.4}, we deduce that
\begin{equation}\begin{aligned}
\|n(\cdot,t)+m(\cdot,t)\|_{L^p(\Omega)}
\leq&\|n(\cdot,t)\|_{L^p(\Omega)}+\|m(\cdot,t)\|_{L^p(\Omega)}\\
\leq&\|n(\cdot,t)\|_{L^p(\Omega)}+K_0|\Omega|^{\frac{1}{p}}
\end{aligned}\label{2.9}\end{equation}
for any $t\in(0,T_{\max}).$ Substituting \dref{2.9} into \dref{2.8} thus shows that
$$\begin{aligned}
\|A^{\alpha}u(\cdot,t)\|_{L^2(\Omega)}\leq& C_1\cdot\left\{1+K_0|\Omega|^{\frac{1}{p}}+\sup_{s\in(0,t)}\|n(\cdot,s)\|_{L^p(\Omega)}\right\}^{\frac{p}{p-1}\cdot\left(\frac{4\alpha-1}{6}+\delta\right)}\\
\leq& C_2\cdot\left\{1+\sup_{s\in(0,t)}\|n(\cdot,s)\|_{L^p(\Omega)}\right\}^{\frac{p}{p-1}\cdot\left(\frac{4\alpha-1}{6}+\delta\right)}
\end{aligned}$$
with $C_2:=C_1\cdot(1+K_0|\Omega|^{\frac{1}{p}})$ for any $t\in(0,T_{\max}),$ which implies \dref{2.6}. According to the arguments of \cite[Corollary 2.1]{WinklerC}, there exists some $C_3=C_3(K_0,|\Omega|,p,q,\delta)>0$ such that
$$\|u(\cdot,t)\|_{L^q(\Omega)}\leq C_3\cdot\left\{1+\sup_{s\in(0,t)}\|n(\cdot,s)+m(\cdot,s)\|_{L^p(\Omega)}\right\}^{\frac{p}{p-1}\cdot\left(\frac{q-3}{3q}+\delta\right)}
$$
for any $t\in(0,T_{\max}),$ from which and \dref{2.9} we can establish \dref{2.7} in a similar manner.

\section{Conditional estimates for the gradients of signal}

The aim of this section is to establish a temporally independent $L^{\infty}$ bound for $\nabla c,$ subject to a conditional estimate appearing as that on the right hand sides of \dref{2.6} and \dref{2.7}. Recalling the arguments pursuing in \cite{WinklerC}, it is essential to resort to higher order conditional estimates as compare to $W^{1,\infty}$-topology, which is based on a combination of the estimates provided by Lemma 2.2 with proper applications of the $L^p$-$L^q$ estimates of the sectorial operator (\cite{Horstmann}). Here and throughout the sequel, we abbreviate $B:=B_l$ to stand for the sectorial operator $-\Delta+1$ under homogeneous Neumann boundary conditions in $\bigcap_{l>1}L^l(\Omega),$ and use $(B^{\mu})_{\mu>0}$ to denote the family of positive fractional powers $B^{\mu}=B^{\mu}_l.$ Moreover, in order to express in more consice forms, we let
\begin{equation}
H_p(t):=1+\sup_{s\in(0,t)}\|n(\cdot,s)\|_{L^p(\Omega)},\quad t\in(0,T_{\max})
\label{3.1}\end{equation}
and
\begin{equation}
J_{l,\mu}(t):=1+\sup_{s\in(0,t)}\left\|B^{\mu}(c(\cdot,s)-e^{-sB}c_0)\right\|_{L^l(\Omega)},\quad t\in(0,T_{\max}).
\label{3.2}\end{equation}

\textbf{Lemma 3.1} \emph{Suppose that $\mu\in(\frac{1}{2},1)$ and $l>3.$  Then for all $\delta>0$ there exists some $C=C(\mu,l,\delta)>0$ such that
\begin{equation}\left\|\nabla(c(\cdot,t)-e^{-sB}c_0)\right\|_{L^{\infty}(\Omega)}\leq C\cdot\left\{1+\sup_{s\in(0,t)}\left\|B^{\mu}(c(\cdot,s)-e^{-sB}c_0)\right\|_{L^l(\Omega)}\right\}^{\frac{l+3}{2\mu l}+\delta}\label{3.3}\end{equation}
for all $t\in(0,T_{\max}).$}

\textbf{Proof.} Since $\mu\in(\frac{1}{2},1)$ allows for choices of $l>3$ sufficiently large and $\delta>0$ arbitrarily small such that $1-\frac{l+3}{2\mu l}>0$ and $\delta<1-\frac{l+3}{2\mu l},$ we can thus take
\begin{equation}
\nu(\delta):=\frac{l+3}{2l}+\delta\mu<\mu.
\label{3.4}\end{equation}
Invoking the interpolation inequality established in \cite[Theorem 2.14.1]{Friedman} for fractional powers of sectorial operators, we can find $C_1=C_1(\mu,l,\delta)>0$ such that
$$\begin{aligned}
\left\|B^{\nu}(c(\cdot,t)-e^{-sB}c_0)\right\|_{L^l(\Omega)}\leq & C_1\left\|B^{\mu}(c(\cdot,t)-e^{-sB}c_0)\right\|^{\frac{\nu}{\mu}}_{L^l(\Omega)}\left\|c(\cdot,t)-e^{-sB}c_0\right\|^{\frac{\mu-\nu}{\mu}}_{L^l(\Omega)}\\
\leq&C_1\left\{2K_0|\Omega|^{\frac{1}{l}}\right\}^{1-\delta-\frac{l+3}{2\mu l}}\left\|B^{\mu}(c(\cdot,t)-e^{-sB}c_0)\right\|^{\frac{l+3}{2\mu l}+\delta}_{L^l(\Omega)}
\end{aligned}$$
with $K_0>0$ as chosen in \dref{1.6} for all $t\in(0,T_{\max}),$ which combined with the continuous embedding $D(B^{\nu})\hookrightarrow W^{1,\infty}(\Omega)$ (\cite{Henry}) provides $C_2=C_2(\mu,l,\delta)>0$ such that
$$\begin{aligned}
\left\|\nabla(c(\cdot,t)-e^{-sB}c_0)\right\|_{L^{\infty}(\Omega)}\leq&
C_2\left\|B^{\nu}(c(\cdot,t)-e^{-sB}c_0)\right\|_{L^l(\Omega)}\\
\leq&C_2C_1\left\{2K_0|\Omega|^{\frac{1}{l}}\right\}^{1-\delta-\frac{l+3}{2\mu l}}\left\|B^{\mu}(c(\cdot,t)-e^{-sB}c_0)\right\|^{\frac{l+3}{2\mu l}+\delta}_{L^l(\Omega)}
\end{aligned}$$
for all $t\in(0,T_{\max}),$ and thus \dref{3.3} follows with $C:=C_2C_1\left\{2K_0|\Omega|^{\frac{1}{l}}\right\}^{1-\delta-\frac{l+3}{2\mu l}}.$

Relying on Lemma 3.1, we can achieve the following conditional estimates in topology of $D(B^{\mu}).$

\textbf{Lemma 3.2} \emph{Let $l>3,p\geq2$ and $\mu\in(\frac{1}{2},1).$  Then for any $\delta>0$ one can find $C=C(l,p,\mu,\delta)>0$ satisfying
\begin{equation}\left\|B^{\mu}(c(\cdot,s)-e^{-sB}c_0)\right\|_{L^l(\Omega)}\leq C\cdot\left\{1+\sup_{s\in(0,t)}\|n(\cdot,s)\|_{L^p(\Omega)}\right\}^{\frac{p}{p-1}\cdot\left(\frac{2\mu}{3}+\delta\right)}
\label{3.5}\end{equation}
for all $t\in(0,T_{\max}).$}

\textbf{Proof.} Taking $\delta>0$ small enough such that
\begin{equation}
\delta<\min\left\{1-\frac{l+3}{2\mu l},2l(1-\mu)\right\},
\label{3.6}\end{equation}
we let
\begin{equation}q:=\frac{3l}{3+2l(1-\mu)-\delta}.\label{3.6a}\end{equation}
Then thanks to $\delta<2l(1-\mu)$ and to $l>3+2l-2l\mu+2l\mu\delta>3+2l-2l\mu$ implied by \dref{3.6}, one can see that
\begin{equation}
l=\frac{3l}{3+2l-2l\mu-2l(1-\mu)}>q=\frac{3l}{3+2l-2l\mu-\delta}>\frac{3l}{3+2l-2l\mu}>3.
\label{3.7}\end{equation}
Now, applying $B^{\mu}$ to both sides of the variation-of-constants representation
$$c(\cdot,t)-e^{-tB}c_0=\int^t_0e^{-B(t-s)}\left\{m(\cdot,s)-u(\cdot,s)\nabla c(\cdot,s)\right\}ds$$
for all $t\in(0,T_{\max}),$ we obtain
$$\left\|B^{\mu}(c(\cdot,t)-e^{-sB}c_0)\right\|_{L^l(\Omega)}
\leq\int^t_0\left\|B^{\mu}e^{-B(t-s)}m(\cdot,s)\right\|_{L^l(\Omega)}ds
+\int^t_0\left\|B^{\mu}e^{-B(t-s)}u(\cdot,s)\nabla c(\cdot,s)\right\|_{L^l(\Omega)}ds$$
for all $t\in(0,T_{\max}),$ where in accordance with the $L^p$-$L^q$ estimates of the sectorial operator (\cite{Horstmann}), it follows from \dref{1.6}, \dref{2.5}, \dref{2.7}, \dref{3.1}, \dref{3.2} and \dref{3.3} that
$$\begin{aligned}
\int^t_0\left\|B^{\mu}e^{-B(t-s)}m(\cdot,s)\right\|_{L^l(\Omega)}ds
\leq&C_3\int^t_0\left(1+(t-s)^{-\mu}\right)e^{-(t-s)}\|m(\cdot,s)\|_{L^l(\Omega)}ds\\
\leq&C_3K_0|\Omega|^{\frac{1}{l}}\int^t_0\left(1+(t-s)^{-\mu}\right)e^{-(t-s)}ds\leq C_4
\end{aligned}$$
for all $t\in(0,T_{\max})$ with some $C_3>0$ and $C_4:=C_3K_0|\Omega|^{\frac{1}{l}}\int^{\infty}_0\left(1+\rho^{-\mu}\right)e^{-\rho}d\rho<\infty$ due to $\mu\in(\frac{1}{2},1),$ and that
$$\begin{aligned}
&\int^t_0\left\|B^{\mu}e^{-B(t-s)}u(\cdot,s)\nabla c(\cdot,s)\right\|_{L^l(\Omega)}ds\\
\leq&C_5\int^t_0\left(1+(t-s)^{-\mu-\frac{3}{2}(\frac{1}{q}-\frac{1}{l})}\right)e^{-(t-s)}\|u(\cdot,s)\nabla c(\cdot,s)\|_{L^q(\Omega)}ds\\
\leq&C_5\int^t_0\left(1+(t-s)^{-\mu-\frac{3}{2}(\frac{1}{q}-\frac{1}{l})}\right)e^{-(t-s)}\|u(\cdot,s)\|_{L^q(\Omega)}
\|\nabla c(\cdot,s)\|_{L^{\infty}(\Omega)}ds\\
\leq&C_5\int^t_0\left(1+(t-s)^{-\mu-\frac{3}{2}(\frac{1}{q}-\frac{1}{l})}\right)e^{-(t-s)}\|u(\cdot,s)\|_{L^q(\Omega)}
\cdot\left\{\left\|\nabla (c(\cdot,s)-e^{-tB}c_0)\right\|_{L^{\infty}(\Omega)}\right.\\
&\left.+\|\nabla e^{-tB}c_0\|_{L^{\infty}(\Omega)}\right\}ds\\
\leq&C_5K_2\int^t_0\left(1+(t-s)^{-\mu-\frac{3}{2}(\frac{1}{q}-\frac{1}{l})}\right)e^{-(t-s)}ds\cdot H^{\frac{p}{p-1}\cdot\left(\frac{q-3}{3q}+\delta\right)}_p(t)\cdot\left\{C_6J^{\frac{l+3}{2\mu l}+\delta}_{l,\mu}(t)+C_7\|\nabla c_0\|_{L^{\infty}(\Omega)}\right\}\\
\leq&C_8H^{\frac{p}{p-1}\cdot\left(\frac{q-3}{3q}+\delta\right)}_p(t)\cdot J^{\frac{l+3}{2\mu l}+\delta}_{l,\mu}(t)
\end{aligned}$$
for all $t\in(0,T_{\max}),$ where $C_8:=C_5K_2(C_6+C_7K_0)\int^{\infty}_0\left(1+{\rho}^{-\mu-\frac{3}{2}(\frac{1}{q}-\frac{1}{l})}\right)e^{-\rho}d\rho<\infty$ thanks to \dref{3.7}. Thereupon,
\begin{equation}
\left\|B^{\mu}(c(\cdot,t)-e^{-sB}c_0)\right\|_{L^l(\Omega)}
\leq C_4+C_8H^{\frac{p}{p-1}\cdot\left(\frac{q-3}{3q}+\delta\right)}_p(t)\cdot J^{\frac{l+3}{2\mu l}+\delta}_{l,\mu}(t)
\label{3.8}\end{equation}
for all $t\in(0,T_{\max}).$ Observing from \dref{3.6} that
$$\delta+\frac{l+3}{2\mu l}<1,$$
we make use of Young's inequality to attain
$$\left\|B^{\mu}(c(\cdot,t)-e^{-sB}c_0)\right\|_{L^l(\Omega)}
\leq C_4+\frac{1}{2}J_{l,\mu}(t)+C_9H^{\frac{p}{p-1}\cdot\left(\frac{q-3}{3q}+\delta\right)\cdot\frac{2\mu l}{2\mu l-l-3-2\mu l\delta}}_p(t)$$
with certain $C_9=C_9(\mu,l,p,q,\delta)$ for all $t\in(0,T_{\max}),$ which further implies
$$J_{l,\mu}(t)\leq 1+C_4+\frac{1}{2}J_{l,\mu}(t)+C_9H^{\frac{p}{p-1}\cdot\left(\frac{q-3}{3q}+\delta\right)\cdot\frac{2\mu l}{2\mu l-l-3-2\mu l\delta}}_p(t),$$
that is
\begin{equation}
J_{l,\mu}(t)\leq 2(1+C_4)+2C_9H^{\frac{p}{p-1}\cdot\left(\frac{q-3}{3q}+\delta\right)\cdot\frac{2\mu l}{2\mu l-l-3-2\mu l\delta}}_p(t)
\label{3.9}\end{equation}
for all $t\in(0,T_{\max}).$ Letting
$$\phi(\tilde\delta):=\frac{p}{p-1}\cdot\left(\frac{2\mu l-l-3+\tilde\delta}{3l}+\tilde\delta\right)\cdot\frac{2\mu l}{2\mu l-l-3-2\mu l\tilde\delta},$$
we find that $\phi(\tilde\delta)\searrow\frac{p}{p-1}\cdot\frac{2\mu}{3}$ as $\tilde\delta\searrow 0,$ whence for some chosen $\delta>0$ it is possible to pick $\delta'\in\left(0,\min\left\{1-\frac{l+3}{2\mu l},2l(1-\mu)\right\}\right)$ such that
$$\phi(\delta')\leq\frac{p}{p-1}\cdot\frac{2\mu}{3}+\delta.$$
Thus, upon an elementary calculation, we draw on \dref{3.6a} to have
$$\begin{aligned}\frac{p}{p-1}\cdot\left(\frac{q-3}{3q}+\delta'\right)\cdot\frac{2\mu l}{2\mu l-l-3-2\mu l\delta'}&=\frac{p}{p-1}\cdot\left(\frac{2\mu l-l-3+\delta'}{3l}+\delta'\right)\cdot\frac{2\mu l}{2\mu l-l-3-2\mu l\delta'}\\
&=\phi(\delta')\leq\frac{p}{p-1}\cdot\frac{2\mu}{3}+\delta,\end{aligned}$$
which along with \dref{3.9} and \dref{3.1} implies \dref{3.5}.

Based on a well-known continuous embedding, a combination of Lemma 3.1 with Lemma 3.2 provides the desired uniform $L^{\infty}$ conditional estimates of $\nabla c$ as follows.

\textbf{Lemma 3.3} \emph{Assume that $p\geq2$ and $\delta>0.$  Then there exists $C(p,\delta)>0$ such that
\begin{equation}\left\|\nabla c(\cdot,t)\right\|_{L^{\infty}(\Omega)}\leq C\cdot\left\{1+\sup_{s\in(0,t)}\|n(\cdot,s)\|_{L^p(\Omega)}\right\}^{\frac{p}{p-1}\cdot\left(\frac{1}{3}+\delta\right)}
\label{3.11}\end{equation}
for all $t\in(0,T_{\max}).$}

\textbf{Proof.} For given $\delta>0,$ we can find $l>3$ large enough such that
$$\frac{1}{l}<\delta,$$
and whereby
$$\frac{l+3}{3l}<\frac{1}{3}+\delta.$$
Let
$$\varphi(\tilde\delta):=\left(\frac{l+3}{2\mu l}+\tilde\delta\right)\cdot\left(\frac{2\mu}{3}+\tilde\delta\right),\qquad\tilde\delta>0.$$
It is evident that
$$\varphi(\tilde\delta)\searrow\frac{l+3}{2\mu l}\cdot\frac{2\mu}{3}=\frac{l+3}{3l}<\frac{1}{3}+\delta\quad\textrm{as}\quad\tilde\delta\searrow 0,$$
which allows for a choice of $\delta''=\delta''(\delta)>0$ such that
\begin{equation}
\varphi(\delta'')\leq\frac{1}{3}+\delta.
\label{3.12}\end{equation}
Combining Lemma 3.1 with Lemma 3.2 provides some $C_1=C_1(p,l,\mu,\delta'')>0$ fulfilling
\begin{equation}
\left\|\nabla(c(\cdot,t)-e^{-sB}c_0)\right\|_{L^{\infty}(\Omega)}\leq C_1H^{\frac{p}{p-1}\cdot\left(\frac{2\mu}{3}+\delta''\right)\cdot\left(\frac{l+3}{2\mu l}+\delta''\right)}_p
\label{3.13}\end{equation}
for all $t\in(0,T_{\max}).$ Moreover, from the regularity properties of the Neumman heat semigroup (\cite{Henry,Winkler 1}), we have
\begin{equation}
\left\|\nabla e^{-tB}c_0\right\|_{L^{\infty}(\Omega)}\leq C_2\left\|\nabla c_0\right\|_{L^{\infty}(\Omega)}
\label{3.14}\end{equation}
for all $t\in(0,T_{\max}).$ Therefore, a collection of \dref{3.1}, \dref{3.12}, \dref{3.13} and \dref{3.14} entails
$$\begin{aligned}\left\|\nabla c(\cdot,t)\right\|_{L^{\infty}(\Omega)}\leq &\left\|\nabla(c(\cdot,t)-e^{-sB}c_0)\right\|_{L^{\infty}(\Omega)}+\left\|\nabla e^{-sB}c_0\right\|_{L^{\infty}(\Omega)}\\
\leq&C_1H^{\frac{p}{p-1}\cdot\left(\frac{2\mu}{3}+\delta''\right)\cdot\left(\frac{l+3}{2\mu l}+\delta''\right)}_p(t)+C_2\left\|\nabla c_0\right\|_{L^{\infty}(\Omega)}\\
\leq&C_3H^{\frac{p}{p-1}\cdot\left(\frac{2\mu}{3}+\delta''\right)\cdot\left(\frac{l+3}{2\mu l}+\delta''\right)}_p(t)\\
=&C_3H^{\frac{p}{p-1}\cdot\varphi(\delta'')}_p(t)\\
\leq&C_3H^{\frac{p}{p-1}\cdot\left(\frac{1}{3}+\delta\right)}_p(t)\end{aligned}$$
with $C_3:=C_1+C_2\left\|\nabla c_0\right\|_{L^{\infty}(\Omega)}$ for all $t\in(0,T_{\max}),$ which shows \dref{3.11}.


\section{Boundedness in \dref{1.4}. Proof of Theorem 1.1}

With the aids of Lemma 3.3, a standard testing procedure as used in the arguments of \cite[Lemma 3.2]{WinklerC} can yield the $L^p$ bounds for the component $n.$

\textbf{Lemma 4.1} \emph{Let both \dref{1.5} and \dref{1.6} hold. If \dref{1.7} is fulfilled with some $K_{\mathcal{S}}>0$ and $\theta>0,$ then one can find some $C=C(p)>0$ such that
\begin{equation}
\sup_{t\in(0,T_{\max})}\|n(\cdot,t)\|_{L^p(\Omega)}\leq C\qquad\textrm{for any}~~p>3.
\label{4.1}\end{equation}}

\textbf{Proof.} Without loss of generality, we suppose $\theta<1,$ then it is possible to take $\delta>0$ sufficiently small such that
\begin{equation}
(1+3\delta)(1-\theta)<1.
\label{4.2}\end{equation}
Testing $n$-equation in \dref{1.4} by $n^{p-1}$ and integrating by parts, we obtain from $\nabla\cdot u=0,$ Young's inequality and the nonnegativity of both $n$ and $m$ that
$$\begin{aligned}
\frac{1}{p}\frac{d}{dt}\int_{\Omega}n^p+(p-1)\int_{\Omega}n^{p-2}|\nabla n|^2=&(p-1)\int_{\Omega}n^{p-1}\mathcal{S}(|\nabla c|^2)\nabla c\cdot\nabla n-\int_{\Omega}n^p m\\
\leq&(p-1)\int_{\Omega}n^{p-1}\mathcal{S}(|\nabla c|^2)\nabla c\cdot\nabla n\\
\leq&\frac{p-1}{2}\int_{\Omega}n^{p-2}|\nabla n|^2+\frac{p-1}{2}\int_{\Omega}n^p\mathcal{S}^2(|\nabla c|^2)|\nabla c|^2
\end{aligned}$$
for all $t\in(0,T_{\max}),$ which is actually
\begin{equation}
\frac{d}{dt}\int_{\Omega}n^p+\frac{2(p-1)}{p}\int_{\Omega}\left|\nabla n^{\frac{p}{2}}\right|^2\leq\frac{p(p-1)K^2_{\mathcal{S}}}{2}\int_{\Omega}n^p|\nabla c|^{2-2\theta}
\label{4.3}\end{equation}
for all $t\in(0,T_{\max}).$ Thanks to $2(1-\theta)>0,$ we draw on Lemma 3.2 to gain
\begin{equation}\begin{aligned}
\frac{p(p-1)K^2_{\mathcal{S}}}{2}\int_{\Omega}n^p|\nabla c|^{2-2\theta}\leq&\frac{p(p-1)K^2_{\mathcal{S}}}{2}\|\nabla c\|^{2-2\theta}_{L^{\infty}(\Omega)}\int_{\Omega}n^p\\
\leq&C_1H^{\frac{p}{p-1}\cdot\left(\frac{1}{3}+\delta\right)\cdot(2-2\theta)}_p(t)\cdot\int_{\Omega}n^p
\end{aligned}\label{4.4}\end{equation}
for all $t\in(0,T_{\max}),$ where $C_1=C_1(p)>0.$ Along with \dref{2.3}, an application of the Gagliardo--Nirenberg inequality provides $C_2=C_2(p)>0$ and $C_3=C_3(p)>0$ such that
\begin{equation}\begin{aligned}
\left\{\int_{\Omega}n^p\right\}^{\frac{3p-1}{3(p-1)}}
=&\|n^{\frac{p}{2}}\|^{\frac{2(3p-1)}{3(p-1)}}_{L^2(\Omega)}\\
\leq&C_2\|\nabla n^{\frac{p}{2}}\|^2_{L^2(\Omega)}\|n^{\frac{p}{2}}\|^{\frac{4}{3(p-1)}}_{L^{\frac{2}{p}}(\Omega)}
+C_2\|n^{\frac{p}{2}}\|^{\frac{2(3p-1)}{3(p-1)}}_{L^{\frac{2}{p}}(\Omega)}\\
\leq&C_3\|\nabla n^{\frac{p}{2}}\|^2_{L^2(\Omega)}+C_3
\end{aligned}\label{4.5}\end{equation}
for all $t\in(0,T_{\max}).$ With
\begin{equation}\lambda:=\frac{3p-1}{3(p-1)},\label{4.5a}\end{equation}
\dref{4.5} implies
\begin{equation}
C_4\cdot\left\{\int_{\Omega}n^p\right\}^{\lambda}-\frac{2(p-1)}{p}\leq\frac{2(p-1)}{p}\int_{\Omega}\left|\nabla n^{\frac{p}{2}}\right|^2
\label{4.6}\end{equation}
for all $t\in(0,T_{\max}),$ where $C_4:=\frac{2(p-1)}{pC_3}.$ Inserting \dref{4.6} and \dref{4.4} into \dref{4.3} and employing Young's inequality yield
$$\begin{aligned}
\frac{d}{dt}\int_{\Omega}n^p+C_4\cdot\left\{\int_{\Omega}n^p\right\}^{\lambda}
\leq&C_1H^{\frac{p}{p-1}\cdot\left(\frac{1}{3}+\delta\right)\cdot(2-2\theta)}_p(T)\cdot\int_{\Omega}n^p+\frac{2(p-1)}{p}\\
\leq&\frac{C_4}{2}\cdot\left\{\int_{\Omega}n^p\right\}^{\lambda}+\left(\frac{2}{C_4}\right)^{\frac{1}{\lambda-1}}\cdot C^{\frac{\lambda}{\lambda-1}}_1\cdot H^{\frac{\lambda}{\lambda-1}\cdot\frac{p}{p-1}\cdot\left(\frac{1}{3}+\delta\right)\cdot(2-2\theta)}_p(T)\\
&+\frac{2(p-1)}{p}\\
\leq&\frac{C_4}{2}\cdot\left\{\int_{\Omega}n^p\right\}^{\lambda}+C_5\cdot H^{\frac{\lambda}{\lambda-1}\cdot\frac{p}{p-1}\cdot\left(\frac{1}{3}+\delta\right)\cdot(2-2\theta)}_p(T)
\end{aligned}\label{4.6}$$
for each $T\in(0,T_{\max})$ with $C_5:=\left(\frac{2}{C_4}\right)^{\frac{1}{\lambda-1}}\cdot C^{\frac{\lambda}{\lambda-1}}_1+\frac{2(p-1)}{p}$ due to $H_p\geq 1,$ and thus
$$\frac{d}{dt}\int_{\Omega}n^p+\frac{C_4}{2}\cdot\left\{\int_{\Omega}n^p\right\}^{\lambda}
\leq C_5\cdot H^{\frac{\lambda}{\lambda-1}\cdot\frac{p}{p-1}\cdot\left(\frac{1}{3}+\delta\right)\cdot(2-2\theta)}_p(T)$$
for all $t\in(0,T).$ By means of an ODE comparison argument, this further entails
\begin{equation}
\int_{\Omega}n^p(\cdot,t)\leq\max\left\{\int_{\Omega}n^p_0,\left\{\frac{2C_5}{C_4}\cdot H^{\frac{\lambda}{\lambda-1}\cdot\frac{p}{p-1}\cdot\left(\frac{1}{3}+\delta\right)\cdot(2-2\theta)}_p(T)\right\}^{\frac{1}{\lambda}}\right\}
\label{4.7}\end{equation}
for all $t\in(0,T).$ In view of \dref{3.1}, we infer from \dref{4.7} that
\begin{equation}\begin{aligned}
H_p(T)\leq& 1+\max\left\{\|n_0\|_{L^p(\Omega)},\left(\frac{2C_5}{C_4}\right)^{\frac{1}{p\lambda}}\cdot H^{\varrho}_p(T)\right\}\\
\leq& C_6 H^{\varrho}_p(T)
\end{aligned}\label{4.8}\end{equation}
for any $T\in(0,T_{\max}),$ where $\varrho:=\frac{1}{\lambda-1}\cdot\frac{1}{p-1}\cdot\left(\frac{1}{3}+\delta\right)\cdot(2-2\theta)$ and $C_6:=1+\max\left\{\|n_0\|_{L^p(\Omega)},\left(\frac{2C_5}{C_4}\right)^{\frac{1}{p\lambda}}\right\}.$
It is clear from \dref{4.5a} and \dref{4.2} that
$$\varrho=\frac{3(p-1)}{2}\cdot\frac{1}{p-1}\cdot\left(\frac{1}{3}+\delta\right)\cdot(2-2\theta)
=(1+3\delta)(1-\theta)<1,$$
whence \dref{4.8} shows that
$$H_p(T)\leq C^{\frac{1}{1-\varrho}}_6\quad\textrm{for any}\quad T\in(0,T_{\max}),$$
which implies \dref{4.1} by letting $T\nearrow T_{\max}.$

Now, we are in the position to pursue the boundedness of each quantity on the left hand side of \dref{2.2}. In particular, the quantities associated with the components of both $c$ and $u$ can be estimated straightforwardly by a collection of Lemmas 2.2--2.3, Lemma 3.3 and Lemma 4.1. For $m,$ the corresponding boundedness needs to be verified by similar strategies as performed in Lemmas 3.2--3.3. In the final, aided by the bounds established for the quantities related to $c,m$ and $u,$ the temporally independent $L^{\infty}$ bounds of $n$ can be achieved through an appropriate application of heat semigroup theories as done in \cite{WinklerC}.

\textbf{Lemma 4.2} \emph{Suppose both \dref{1.5} and \dref{1.6} are valid. Let \dref{1.7} hold with $K_{\mathcal{S}}>0$ and $\theta>0.$ Then there exists $C>0$ with the properties that
\begin{equation}
\|c(\cdot,t)\|_{W^{1,\infty}(\Omega)}\leq C\quad\textrm{for any}\quad t\in(0,T_{\max})
\label{4.10}\end{equation}
and
\begin{equation}
\|m(\cdot,t)\|_{W^{1,\infty}(\Omega)}\leq C\quad\textrm{for any}\quad t\in(0,T_{\max})
\label{4.11}\end{equation}
as well as
\begin{equation}
\|A^{\alpha}u(\cdot,t)\|_{L^2(\Omega)}\leq C\quad\textrm{for any}\quad t\in(0,T_{\max})
\label{4.12}\end{equation}
with certain $\alpha\in(\frac{3}{4},1).$}

\textbf{Proof.} Thanks to Lemma 4.1, both \dref{4.10} and \dref{4.12} are immediate consequences of Lemmas 2.2--2.3 and Lemma 3.3. Specially, due to $\alpha\in(\frac{3}{4},1),$ the continuous embedding together with \dref{4.12} provides $C_1>0$ and $C_2>0$ such that
\begin{equation}
\|u(\cdot,t)\|_{L^{\infty}(\Omega)}\leq C_1\|A^{\alpha}u(\cdot,t)\|_{L^2(\Omega)}\leq C_2
\label{4.13}\end{equation}
for any $t\in(0,T_{\max}).$ For each given $r>3,$ we let $\beta\in(\frac{1}{2},1)$ fulfill
$$\beta>\frac{r+3}{2r},$$
then there exists $\vartheta\in(\frac{1}{2},1)$ satisfying
\begin{equation}
\beta>\vartheta>\frac{r+3}{2r}.
\label{4.14}\end{equation}
Now, applying $B^{\beta}$ on both sides of the Duhamel representation of $m$ and invoking the $L^p$-$L^q$ estimates of the sectorial operator (\cite{Horstmann}), we obtain $C_3>0$ such that
\begin{equation}\begin{aligned}
\left\|B^{\beta}\left(m(\cdot,t)-e^{-Bt}m_0\right)\right\|_{L^r(\Omega)}
\leq&\int^t_0\left\|B^{\beta}e^{-B(t-s)}u(\cdot,s)\nabla m(\cdot,s)\right\|_{L^r(\Omega)}ds\\
&+\int^t_0\left\|B^{\beta}e^{-B(t-s)}m(\cdot,s)(1-n(\cdot,s))\right\|_{L^r(\Omega)}ds\\
\leq&C_3\int^t_0\left(1+(t-s)^{-\beta}\right)e^{-(t-s)}\left\|u(\cdot,s)\nabla m(\cdot,s)\right\|_{L^r(\Omega)}ds\\
&+C_3\int^t_0\left(1+(t-s)^{-\beta}\right)e^{-(t-s)}\left\|m(\cdot,s)(1-n(\cdot,s))\right\|_{L^r(\Omega)}ds
\end{aligned}\label{4.15}\end{equation}
for all $t\in(0,T_{\max}).$ Since \dref{4.14} ensures the inequality
$$\left\|\nabla\left(m(\cdot,t)-e^{-Bt}m_0\right)\right\|_{L^r(\Omega)}\leq C_4\left\|B^{\vartheta}\left((m(\cdot,t)-e^{-Bt}m_0\right)\right\|_{L^r(\Omega)}$$
with some $C_4>0$ due to the embedding $D(B^{\beta})\hookrightarrow W^{1,r}(\Omega)$ for any $t\in(0,T_{\max}),$ and also allows for an application of the following interpolation features of the fractional power of sectorial operators (\cite[Theorem 2.14.1]{Friedman})
$$\left\|B^{\vartheta}\left(m(\cdot,t)-e^{-Bt}m_0\right)\right\|_{L^r(\Omega)}
\leq\left\|B^{\beta}\left(m(\cdot,t)-e^{-Bt}m_0\right)\right\|^{\frac{\vartheta}{\beta}}_{L^r(\Omega)}
\left\|m(\cdot,t)-e^{-Bt}m_0\right\|^{\frac{\beta-\vartheta}{\beta}}_{L^r(\Omega)}$$
for all $t\in(0,T_{\max}),$ we recall \dref{1.6}, \dref{2.4}, \dref{4.14} and make use of the regularity properties of $(e^{-tB})_{t\geq0}$ (\cite{Winkler 1}) to have
\begin{equation}\begin{aligned}
\left\|u(\cdot,t)\nabla m(\cdot,t)\right\|_{L^r(\Omega)}\leq&\left\|u(\cdot,t)\right\|_{L^{\infty}(\Omega)}
\cdot\left\{\left\|\nabla\left(m(\cdot,t)-e^{-Bt}m_0\right)\right\|_{L^r(\Omega)}+\left\|\nabla e^{-Bt}m_0\right\|_{L^r(\Omega)}\right\}\\
\leq&C_2C_4\left\|B^{\vartheta}(m(\cdot,t)-e^{-Bt}m_0)\right\|_{L^r(\Omega)}+C_2C_5\left\|\nabla m_0\right\|_{L^r(\Omega)}\\
\leq&C_2C_4\left\|B^{\beta}\left(m(\cdot,t)-e^{-Bt}m_0\right)\right\|^{\frac{\vartheta}{\beta}}_{L^r(\Omega)} \left(\left\|m(\cdot,t)\right\|_{L^r(\Omega)}+\left\|m_0\right\|_{L^r(\Omega)}\right)^{\frac{\beta-\vartheta}{\beta}}\\ &+C_2C_5K_0|\Omega|^{\frac{1}{r}}\\
\leq&C_2C_4\left(2K_0|\Omega|^{\frac{1}{r}}\right)^{\frac{\beta-\vartheta}{\beta}}\left\|B^{\beta}\left(m(\cdot,t)-e^{-Bt}m_0\right)\right\|^{\frac{\vartheta}{\beta}}_{L^r(\Omega)} +C_2C_5K_0|\Omega|^{\frac{1}{r}}\\
\leq&C_6\left\|B^{\beta}\left(m(\cdot,t)-e^{-Bt}m_0\right)\right\|^{\frac{\vartheta}{\beta}}_{L^r(\Omega)}+C_6
\end{aligned}\label{4.16}\end{equation}
with some $C_5>0$ and $C_6:=\max\left\{C_2C_4\left(2K_0|\Omega|^{\frac{1}{r}}\right)^{\frac{\beta-\vartheta}{\beta}},C_2C_5K_0|\Omega|^{\frac{1}{r}}\right\}$ for all $t\in(0,T_{\max}).$ Moreover, Lemma 4.1 along with \dref{2.4} shows
\begin{equation}\begin{aligned}
\left\|m(\cdot,t)(1-n(\cdot,t))\right\|_{L^r(\Omega)}
\leq&\left\|m(\cdot,t)\right\|_{L^{\infty}(\Omega)}\cdot\left(|\Omega|^{\frac{1}{r}}+\left\|n(\cdot,t)\right\|_{L^r(\Omega)}\right)\\
\leq&\left\|m_0\right\|_{L^{\infty}(\Omega)}\cdot\left(|\Omega|^{\frac{1}{r}}+C_7\right)\\
\leq&K_0\cdot\left(|\Omega|^{\frac{1}{r}}+C_7\right)\\
=:&C_8
\end{aligned}\label{4.17}\end{equation}
for all $t\in(0,T_{\max}),$ where $C_7:=C_7(r)>0.$ Inserting \dref{4.16} and \dref{4.17} into \dref{4.15} entails
\begin{equation}\begin{aligned}
&\left\|B^{\beta}\left(m(\cdot,t)-e^{-Bt}m_0\right)\right\|_{L^r(\Omega)}\\
\leq&C_3(C_6+C_8)\cdot\int^t_0\left(1+(t-s)^{-\beta}\right)e^{-(t-s)}ds\\
&+C_3C_6\cdot\sup_{\tau\in(0,t)}\left\|B^{\beta}\left(m(\cdot,\tau)-e^{-B\tau}m_0\right)\right\|^{\frac{\vartheta}{\beta}}_{L^r(\Omega)} \cdot\int^t_0\left(1+(t-s)^{-\beta}\right)e^{-(t-s)}ds\\
\leq&C_9+C_9\cdot\sup_{\tau\in(0,t)}\left\|B^{\beta}\left(m(\cdot,\tau)-e^{-B\tau}m_0\right)\right\|^{\frac{\vartheta}{\beta}}_{L^r(\Omega)}
\end{aligned}\label{4.18}\end{equation}
with $C_9:=C_3(C_6+C_8)\int^{\infty}_0\left(1+\rho^{-\beta}\right)e^{-\rho}d\rho<\infty$ thanks to $\beta<1$ for all $t\in(0,T_{\max}).$ Define
$$M_{\beta,r}(t):=1+ \sup_{s\in(0,t)}\left\|B^{\beta}\left(m(\cdot,t)-e^{-Bt}m_0\right)\right\|_{L^r(\Omega)}$$
for all $t\in(0,T_{\max}).$ Then \dref{4.18} implies
$$\begin{aligned}
M_{\beta,r}(t)\leq& 1+2C_9\cdot M^{\frac{\vartheta}{\beta}}_{\beta,r}(t)\\
\leq&C_{10}\cdot M^{\frac{\vartheta}{\beta}}_{\beta,r}(t)
\end{aligned}$$
for all $t\in(0,T_{\max}),$ where $C_{10}:=1+2C_9,$ and thus from \dref{4.14},
\begin{equation}
M_{\beta,r}(t)\leq C^{\frac{\beta}{\beta-\vartheta}}_{10}
\label{4.19}\end{equation}
for all $t\in(0,T_{\max}),$ which combined with the embedding $D(B^{\beta})\hookrightarrow W^{1,\infty}(\Omega)$ provides $C_{11}=C_{11}(r,\beta,K_0,|\Omega|)>0$ such that
$$\left\|\nabla\left(m(\cdot,t)-e^{-Bt}m_0\right)\right\|_{L^{\infty}(\Omega)}\leq C_{11}$$
for all $t\in(0,T_{\max}).$ As a result, by \dref{1.6} and the heat semigroup estimates (\cite{Winkler 1}), we achieve
$$\begin{aligned}\left\|\nabla m(\cdot,t)\right\|_{L^{\infty}(\Omega)}\leq& \left\|\nabla\left(m(\cdot,t)-e^{-Bt}m_0\right)\right\|_{L^{\infty}(\Omega)}+\left\|\nabla e^{-Bt}m_0\right\|_{L^{\infty}(\Omega)}\\
\leq&C_{11}+C_{12}\left\|\nabla m_0\right\|_{L^{\infty}(\Omega)}\\
\leq&C_{11}+C_{12}K_0
\end{aligned}$$
with some $C_{12}>0$ for all $t\in(0,T_{\max}).$ This together with \dref{2.4} yields \dref{4.11}.

\textbf{Lemma 4.3} \emph{If both \dref{1.5} and \dref{1.6} are satisfied, and if \dref{1.7} is fulfilled with some $K_{\mathcal{S}}>0$ and $\theta>0,$ then one can find $C>0$ such that
\begin{equation}
\|n(\cdot,t)\|_{L^{\infty}(\Omega)}\leq C\quad\textrm{for each}\quad t\in(0,T_{\max}).
\label{4.20}\end{equation}}

\textbf{Proof.} Let $g_1:=n\mathcal{S}(|\nabla c|^2)\nabla c+un$ and $g_2:=n(1-m).$ Then for each $\iota>3$ Lemmas 4.1--4.2 combined with \dref{2.4} guarantee the existence of $C_1=C_1(\iota)>0$ and $C_2=C_2(\iota)>0$ such that
$$\|g_1(\cdot,t)\|_{L^{\iota}(\Omega)}\leq C_1\quad\textrm{and}\quad \|g_2(\cdot,t)\|_{L^{\frac{\iota}{2}}(\Omega)}\leq C_2\quad\textrm{for each}\quad t\in(0,T_{\max}).$$
Thereupon, proceeding along a similar reasoning as that of \cite[Lemma 3.4]{WinklerC}, we establish \dref{4.20}.

\textbf{Proof of Theorem 1.1.} In view of the blow-up criterion \dref{2.2}, Theorem 1.1 is a direct consequence of Lemmas 4.2--4.3.

\section{Equilibration. Proof of Theorem 1.2}

\subsection{Convergence in $L^2(\Omega)$}

As the cornerstone of this section, the following assertions on decay properties of $nm$ and $\nabla m$ as well as on the stabilization of spatial $L^1$-integrals of both $n$ and $m$ follow from a straightforward testing procedure together with an argument analogous to that of \cite[Lemma 4.2]{Espejo1}.

\textbf{Lemma 5.1} \emph{The components $n$ and $m$ of the solutions constructed in Theorem 1.1 have the properties that
\begin{equation}\int^{\infty}_0\int_{\Omega} n m < \infty,\quad\int^{\infty}_0\int_{\Omega} |\nabla m|^2 < \infty,\label{5.1}\end{equation}
and
\begin{equation}\int_{\Omega} n (\cdot,t)\rightarrow\left\{\int_{\Omega} n_0-\int_{\Omega} m_0\right\}_{+},\quad\int_{\Omega} m (\cdot,t)\rightarrow\left\{\int_{\Omega} m_0-\int_{\Omega} n_0\right\}_{+}\label{5.2}\end{equation}
as $t\rightarrow\infty.$}

\textbf{Proof.} Integrating $n$-equation and $m$-equation over $\Omega\times(0,t)$ for any $t>0,$ respectively, we obtain from $n\geq 0$ and $m\geq 0$ that
$$\int^t_0\int_{\Omega} n m\leq \min\left\{\int_{\Omega} n_0,\int_{\Omega} m_0\right\},$$
which implies the first inequality in \dref{5.1}. Testing $m$-equation by $m$ and integrating the resulted equation on $(0,t)$ yield
$$\frac{1}{2}\int_{\Omega}m^2 (\cdot,t)+\int^t_0\int_{\Omega} |\nabla m|^2=\frac{1}{2}\int_{\Omega}m^2_0-\int^t_0\int_{\Omega} nm^2\quad\textrm{for any}\quad t>0,$$
whence the second inequality in \dref{5.1} holds thanks to the nonnegativity of $n$ and $m.$ Relying on \dref{5.1}, the convergence in \dref{5.2} can be obtained in accordance with the reasoning of \cite[Lemma 4.2]{Espejo1}.

Besides the decay features and the convergence involved in Lemma 5.1, some higher order estimates, such as the H\"{o}lder estimates, are also essential for the derivation of the desired convergence.

\textbf{Lemma 5.2} \emph{For components $n$ and $m,$ one can find some $\gamma\in(0,1)$ and $C>0$ such that
\begin{equation}\|n\|_{C^{\gamma,\frac{\gamma}{2}}(\bar\Omega\times[t,t+1])}\quad\textrm{for all}\quad t>1\label{5.3}\end{equation}
and
\begin{equation}\|m\|_{C^{\gamma,\frac{\gamma}{2}}(\bar\Omega\times[t,t+1])}\quad\textrm{for all}\quad t>1.\label{5.4}\end{equation}}

\textbf{Proof.} Denoting $\xi_1:=n\mathcal{S}(|\nabla c|^2)\nabla c+un$ and $\xi_2:=nm,$ $n$-equation in \dref{1.4} can be rewritten as
$$n_t=\Delta n-\nabla\cdot\xi_1(x,t)-\xi_2(x,t),\qquad x\in\Omega,~~t>0,$$
thanks to $\nabla\cdot u=0,$ where from \dref{4.20}, \dref{4.10}, \dref{1.7}, \dref{4.13} and \dref{2.5}, we can infer that both $\xi_1$ and $\xi_2$ are bounded in $\Omega\times(0,\infty),$ which implies \dref{5.3} according to \cite[Theeorem 1.3]{Porzio}. Similarly, we let $\zeta:=mu,$ then again by $\nabla\cdot u=0$ another equivalent expression of $m$-equation appears as
$$m_t=\Delta m-\nabla\cdot\zeta(x,t)-\xi_2(x,t),\qquad x\in\Omega,~~t>0.$$
It is evident from \dref{2.4} and \dref{4.13} that $\zeta$ is bounded in $\Omega\times(0,\infty),$ whereupon combining with the boundedness feature of $\xi_2$ we conclude that \dref{5.4} also holds.

In light of Lemmas 5.1--5.2, let us provide the uniform $L^2(\Omega)$-convergence of $m$ at first.

\textbf{Lemma 5.3} \emph{For component $m,$ we have
\begin{equation}m(\cdot,t)\rightarrow\frac{1}{|\Omega|}\left\{\int_{\Omega} m_0-\int_{\Omega} n_0\right\}_{+}\quad\textrm{in}~L^2(\Omega)~\textrm{as}~ t\rightarrow\infty.\label{5.5}\end{equation}}

\textbf{Proof.} From \dref{5.1} and the Poincar\'{e} inequality, we have
\begin{equation}
\int^{\infty}_0\left\|m(\cdot,t)-\overline{m}(\cdot,t)\right\|^2_{L^2(\Omega)}dt\leq C_1\int^{\infty}_0\int_{\Omega}|\nabla m|^2<\infty,
\label{5.6}\end{equation}
where $C_1>0$ and $\overline{m}(\cdot,t):=\frac{1}{|\Omega|}\int_{\Omega} m(\cdot,t)$ for $t\geq 0.$ Moreover, the H\"{o}lder continuity of $m$ in \dref{5.4} implies the uniform continuity of $0\leq t\rightarrow\left\|m(\cdot,t)-\overline{m}(\cdot,t)\right\|^2_{L^2(\Omega)}.$ Thus, by means of a reasoning similar to that of \cite[Theorem 1.1]{Tomasz}, we conclude from \dref{5.6} that
$$\left\|m(\cdot,t)-\overline{m}(\cdot,t)\right\|^2_{L^2(\Omega)}\rightarrow 0\qquad\textrm{as}~ t\rightarrow\infty,$$
whence there exists some $t_1>0$ such that for arbitrary $\varepsilon>0$
\begin{equation}
\left\|m(\cdot,t)-\overline{m}(\cdot,t)\right\|^2_{L^2(\Omega)}<\frac{\varepsilon}{4}\qquad\textrm{for all}~ t>t_1.
\label{5.7}\end{equation}
Apart from that, the convergence in \dref{5.2} enable us to choose some $t_2>0$ with the property that
\begin{equation}
\left|\overline{m}(\cdot,t)-\frac{1}{|\Omega|}\left\{\int_{\Omega} m_0-\int_{\Omega} n_0\right\}_{+}\right|^2\leq\frac{\varepsilon}{4|\Omega|}\qquad\textrm{for all}~ t>t_2.
\label{5.8}\end{equation}
Setting $t_0:=\max\{t_1,t_2\}$ and $m_{\infty}:=\frac{1}{|\Omega|}\left\{\int_{\Omega} m_0-\int_{\Omega} n_0\right\}_{+},$ we obtain from \dref{5.7} and \dref{5.8} that
$$\begin{aligned}\int_{\Omega}\left|m(\cdot,t)-m_{\infty}\right|^2
\leq& 2\int_{\Omega}\left|m(\cdot,t)-\overline{m}(\cdot,t)\right|^2+2\int_{\Omega}\left|\overline{m}(\cdot,t)- m_{\infty}\right|^2\\
\leq&\frac{\varepsilon}{2}+\frac{\varepsilon}{2}=\varepsilon\qquad\textrm{for all}~ t>t_0,\end{aligned}$$
which implies \dref{5.5}.

Now, the corresponding convergence of $n,m$ and $u$ in $L^2(\Omega)$ can be achieved by virtue of standard testing procedures as used in \cite{Espejo1}.

\textbf{Lemma 5.4} \emph{For $n,c$ and $u,$ we have
\begin{equation}n(\cdot,t)\rightarrow\frac{1}{|\Omega|}\left\{\int_{\Omega} n_0-\int_{\Omega} m_0\right\}_{+}\quad\textrm{in}~L^2(\Omega)~\textrm{as}~ t\rightarrow\infty\label{5.9}\end{equation}
and
\begin{equation}c(\cdot,t)\rightarrow\frac{1}{|\Omega|}\left\{\int_{\Omega} m_0-\int_{\Omega} n_0\right\}_{+}\quad\textrm{in}~L^2(\Omega)~\textrm{as}~ t\rightarrow\infty\label{5.10}\end{equation}
as well as
\begin{equation}u(\cdot,t)\rightarrow 0\quad\textrm{in}~L^2(\Omega)~\textrm{as}~ t\rightarrow\infty.\label{5.11}\end{equation}}

\textbf{Proof.} The convergence in \dref{5.10} is a immediate consequence of \cite[Lemma 4.7]{Espejo1}. Let $\bar n(\cdot,t):=\frac{1}{|\Omega|}\int_{\Omega} n(\cdot,t)$ for $t\geq 0.$ Then from $n$-equation and $\nabla\cdot u=0,$ we integrate by parts to obtain
\begin{equation}\begin{aligned}
\frac{1}{2}\frac{d}{dt}\int_{\Omega}\left(n(\cdot,t)-\bar n(\cdot,t)\right)^2=&\int_{\Omega}(n(\cdot,t)-\bar n(\cdot,t))\cdot\left(n_t-(\bar n)_t\right)\\
=&\int_{\Omega}\left(n-\bar n\right)\cdot\left(\Delta n-\nabla\cdot(n\mathcal{S}(|\nabla c|^2)\nabla c)-nm-u\cdot\nabla n+\overline{nm}\right)\\
=&-\int_{\Omega}|\nabla n|^2+\int_{\Omega}n\mathcal{S}(|\nabla c|^2)\nabla c\cdot\nabla n-\int_{\Omega}(n-\bar n)\cdot nm+\int_{\Omega}(n-\bar n)\cdot \overline{nm}
\end{aligned}\label{5.12}\end{equation}
for each $t>0,$ where combined with \dref{1.7} and \dref{4.20} an application of Young's inequality entails
$$\begin{aligned}
\int_{\Omega}n\mathcal{S}(|\nabla c|^2)\nabla c\cdot\nabla n
\leq&\frac{1}{2}\int_{\Omega}|\nabla n|^2+\frac{K^2_{\mathcal{S}}}{2}\int_{\Omega}n^2|\nabla c|^2\\
\leq&\frac{1}{2}\int_{\Omega}|\nabla n|^2+\frac{C_1}{2}\int_{\Omega}|\nabla c|^2
\end{aligned}$$
with some $C_1>0$ for each $t>0.$ Thereafter, \dref{5.9} follows from the arguments of \cite[Lemma 4.8]{Espejo1}. In the final, with the aids of \dref{5.9} and \dref{5.10}, one can see that \dref{5.11} is actually valid according to the reasoning of \cite[Lemma 4.9]{Espejo1}.

\subsection{Convergence in the sense of \dref{1.10}. Proof of Theorem 1.2}

Based on the boundedness properties in spaces with higher order regularity and on the convergence previously achieved in $L^2(\Omega),$ we can make use of an Ehrling type lemma to improve the regularity of the spaces, in which each component of the solution converges to corresponding constant equilibrium, to the level as claimed in Theorem 1.2.

\textbf{Proof of Theorem 1.2} In view of \dref{5.3} and \dref{4.12} with $\alpha\in(\frac{3}{4},1),$ there exist $\gamma\in(0,2\alpha-\frac{3}{2})$ and $C_1>0$ such that
\begin{equation}\|n(\cdot,t)-n_{\infty}\|_{C^{\gamma}(\bar\Omega)}\leq C_1\qquad\textrm{for any}~~t>1\label{5.12a}\end{equation}
and
\begin{equation}\|u(\cdot,t)\|_{C^{\gamma}(\bar\Omega)}\leq C_1\qquad\textrm{for any}~~t>0\label{5.12b}\end{equation}
due to the continuous embedding $D(A^{\alpha})\hookrightarrow C^{\gamma}(\bar\Omega)$ (\cite{Henry}), where $n_{\infty}:=\frac{1}{|\Omega|}\left\{\int_{\Omega} n_0-\int_{\Omega} m_0\right\}_{+}.$ Moreover, from \dref{3.5}, \dref{4.1} and \dref{4.19}, we can infer that for some $r>3$ and fixed $\mu\in(\frac{r+3}{2r},1)$ there exists $C_2=C_2(r)>0$ such that
$$\left\|B^{\mu}\left(c(\cdot,t)-e^{-Bt}c_0\right)\right\|_{L^r(\Omega)}\leq C_2\qquad\textrm{for any}~~t>0$$
and
$$\left\|B^{\mu}\left(m(\cdot,t)-e^{-Bt}m_0\right)\right\|_{L^r(\Omega)}\leq C_2\qquad\textrm{for any}~~t>0,$$
and whereby combining with \dref{1.6}, we can find $C_3>0$ satisfying
\begin{equation}\begin{aligned}
\left\|B^{\mu}c(\cdot,t)\right\|_{L^r(\Omega)}
\leq&\left\|B^{\mu}\left(c(\cdot,t)-e^{-Bt}c_0\right)\right\|_{L^r(\Omega)} +\left\|B^{\mu}e^{-Bt}c_0\right\|_{L^r(\Omega)}\\
\leq&C_2+t^{-\mu}e^{-t}\|c_0\|_{L^r(\Omega)}\\
\leq&C_3\qquad\textrm{for any}~~t>1,
\end{aligned}\label{5.13}\end{equation}
and similarly,
\begin{equation}\begin{aligned}
\left\|B^{\mu}m(\cdot,t)\right\|_{L^r(\Omega)}
\leq&\left\|B^{\mu}\left(m(\cdot,t)-e^{-Bt}m_0\right)\right\|_{L^r(\Omega)} +\left\|B^{\mu}e^{-Bt}m_0\right\|_{L^r(\Omega)}\\
\leq&C_2+t^{-\mu}e^{-t}\|m_0\|_{L^r(\Omega)}\\
\leq&C_3\qquad\textrm{for any}~~t>1.
\end{aligned}\label{5.14}\end{equation}
Since $\mu\in(\frac{r+3}{2r},1)$ admits the continuity of the embedding $D(B^{\mu})\hookrightarrow C^{1+\gamma}(\bar\Omega)$ for $\gamma\in(0,2\mu-\frac{r+3}{r}),$ both \dref{5.13} and \dref{5.14} imply the existence of $C_4>0$ such that
$$\|c(\cdot,t)\|_{C^{1+\gamma}(\bar\Omega)}\leq C_4\quad\text{and}\quad\|m(\cdot,t)\|_{C^{1+\gamma}(\bar\Omega)}\leq C_4\qquad\textrm{for any}~~t>1.$$
Therefore, with $m_{\infty}:=\frac{1}{|\Omega|}\left\{\int_{\Omega} m_0-\int_{\Omega} n_0\right\}_{+},$ we can infer that
\begin{equation}
\|c(\cdot,t)-m_{\infty}\|_{C^{1+\gamma}(\bar\Omega)}\leq C_5\quad\text{and}\quad\|m(\cdot,t)-m_{\infty}\|_{C^{1+\gamma}(\bar\Omega)}\leq C_5
\label{5.15}\end{equation}
with certain $C_5>0$ for any $t>1.$ Observing from $C^{\gamma}(\bar\Omega)\hookrightarrow L^{\infty}(\Omega)\hookrightarrow L^2(\Omega)$ and $C^{1+\gamma}(\bar\Omega)\hookrightarrow W^{1,\infty}(\Omega)\hookrightarrow L^2(\Omega)$ that the first embedding of each is compact, we thus apply an Ehrling lemma to obtain some $C_6>0$ such that for any given $\eta>0$
\begin{equation}
\|w\|_{L^{\infty}(\Omega)}\leq\frac{\eta}{2C_1}\|w\|_{C^{\gamma}(\bar\Omega)}+C_6\|w\|_{L^2(\Omega)}\qquad\textrm{for each}~~w\in C^{\gamma}(\bar\Omega)
\label{5.16}\end{equation}
and
\begin{equation}
\|w\|_{W^{1,\infty}(\Omega)}\leq\frac{\eta}{2C_5}\|w\|_{C^{1+\gamma}(\bar\Omega)}+C_6\|w\|_{L^2(\Omega)}\qquad\textrm{for each}~~w\in C^{1+\gamma}(\bar\Omega).
\label{5.17}\end{equation}
Since Lemmas 5.3--5.4 allow for a choice of $t_0>1$ having the properties that
$$\|n(\cdot,t)-n_{\infty}\|_{L^2(\Omega)}\leq\frac{\eta}{2C_6}\qquad\textrm{for any}~~t>t_0,$$
$$\|c(\cdot,t)-m_{\infty}\|_{L^2(\Omega)}\leq\frac{\eta}{2C_6}\qquad\textrm{for any}~~t>t_0$$
and
$$\|m(\cdot,t)-m_{\infty}\|_{L^2(\Omega)}\leq\frac{\eta}{2C_6}\qquad\textrm{for any}~~t>t_0$$
as well as
$$\|u(\cdot,t)\|_{L^2(\Omega)}\leq\frac{\eta}{2C_6}\qquad\textrm{for any}~~t>t_0,$$
\dref{1.10} thereby follows from \dref{5.12a}, \dref{5.12b} and \dref{5.15} in conjunction with applications of \dref{5.16} to both $w=n-n_{\infty}$ and $w=u$ as well as of \dref{5.17} to both $w=c-m_{\infty}$ and $w=m-m_{\infty}.$

\section*{Acknowledgments}

The author is supported by the National Natural Science Foundation of China (Grant No. 11901298), the Fundamental Research Funds for the Central Universities (Grant No. KJQN202052), and the Basic Research Program of Jiangsu Province (Grant No. BK20190504).


\def\refname{References}





\end{document}